\documentclass[12pt,reqno,a4paper]{amsart}
\usepackage[latin1]{inputenc}
\usepackage[T1]{fontenc}
\usepackage{amssymb}
\usepackage{bm}
\usepackage[alpine]{ifsym}

\usepackage[a4paper,top=3.3cm,bottom=3cm,left=3cm,right=3cm,bindingoffset=5mm]{geometry}

\usepackage{color}

\usepackage[colorlinks=true]{hyperref}

\theoremstyle{definition}

\numberwithin{equation}{section}

\begin{document}
\title[Local and global Zeta functions]{The universal zeta function for curve singularities and its relation with global zeta functions}

\author[J.J. Moyano-Fern\'andez]{Julio Jos\'e Moyano-Fern\'andez}

\address{Universitat Jaume I, Campus de Riu Sec, Departamento de Matem\'aticas \& Ins\-titut Universitari de Matem\`atiques i Aplicacions de Castell\'o, 12071
Caste\-ll\'on de la Plana, Spain}  \email{moyano@uji.es}

\thanks{The author was partially supported by the Spanish Government Ministerio de Econom\'ia, Industria y Competitividad (MINECO), grants MTM2012-36917-C03-03 and MTM2015-65764-C3-2-P, as well as by Universitat Jaume I, grant P1-1B2015-02.}

\begin{abstract}
The purpose of this note is to give a brief overview on zeta functions of curve singularities and to provide some evidences on how these and global zeta functions associated to singular algebraic curves over perfect fields relate to each other.
 \end{abstract}

\maketitle

\section{Introduction}\label{intro}

\subsection{}\label{pa:I1} Let $X$ be a complete, geometrically irreducible, singular algebraic curve defined over a perfect field $k$; from now on we will refer to such a curve simply as `algebraic curve over $k$' . Let $K$ be the field of rational functions on $X$. Extending previous works of V. M. Galkin and B. Green---and based on the classical results of F. K. Schmidt \cite{Sch} for nonsingular curves---K.O. St\"ohr (cf. \cite{St1}, \cite{St2}) managed to attach a zeta function to $X$ for finite $k$ in the following manner: If $\mathcal{O}_X$ is the structure sheaf of $X$, he defined the Dirichlet series
$$
\zeta (\mathcal{O}_X,s) :=\sum_{\mathfrak{a} \succeq \mathcal{O}_X} q^{-s\deg{\mathfrak{a}}}, \ \ \ \ s \in \mathbb{C} \ \mbox{with} \ \mathrm{Re}(s) >0, 
$$
where the sum is taken over all positive divisors of $X$, and $\deg{(~\cdot~)}$ denotes the degree of those divisors. Observe that the change of variables $T=q^{-s}$ allows to consider the formal power series in $T$
$$
Z(\mathcal{O}_X,T)=\sum_{n=0}^{\infty} \# (\{ \mbox{positive divisors of } X \mbox{ of degree } n\})\cdot T^{n}.
$$
Moreover, St\"ohr considered local zeta functions, i.e., zeta functions attached to every local ring $\mathcal{O}_P$ of points $P$ at $X$ of the form
$$
Z(\mathcal{O}_P,T):=\sum_{\mathfrak{a}\subseteq \mathcal{O}_P} T^{\deg{\mathfrak{a}}}=\sum_{n=0}^{\infty} \# (\{ \mathrm{positive~}  \mathcal{O}_P\mathrm{-ideals~of~degree~} n\})\cdot  T^{n}.
$$
%\begin{align*}
%Z(\mathcal{O}_P,t):=&\sum_{\mathfrak{a}\subseteq \mathcal{O}_P} t^{\deg{\mathfrak{a}}}\\
%=&\sum_{n=0}^{\infty}\{ \mbox{positive }  \mathcal{O}-\mbox{ideals of degree } n\} t^{n}.
%\end{align*}
This series extends previous definitions by Galkin \cite{G} and Green \cite{Gr}. Furthermore, the Euler product formula for the formal power series yields the identity 
$$
Z(\mathcal{O}_X,T)=\prod_{P\in X} Z(\mathcal{O}_P,T),
$$
which actually establishes a link between the local and global theory.
Every local factor $Z(\mathcal{O}_P,T)$ splits again into factors 
$$
Z(\mathcal{O}_P, \mathcal{O}_P,T)=\sum_{n=0}^{\infty} \# (\{ \mathrm{principal~integral~}  \mathcal{O}\mathrm{-ideals~of~codimension~} n\})\cdot  T^{n}
$$
which are determined by the value semigroup of $\mathcal{O}_P$ (see \S 2.1 for the definition of this semigroup) if the field is big enough, as Z\'u\~niga showed in \cite{Z}.

\subsection{} On the other hand, when studying the Gorenstein property of one-dimensio\-nal local Cohen-Macaulay rings, Campillo, Delgado and Kiyek \cite[(3.8)]{CDK} observed the existence of a Laurent series---a polynomial in their situation---attached to those rings, and satisfying a functional equation in the case of Gorenstein rings. Further investigations by Campillo, Delgado and Gusein-Zade \cite{CDG1}--\cite{CDG9} led to the definition of a Poincar\'e series associated to a complex curve singularity as an integral with respect to the Euler characteristic (se also O. Viro \cite{V}). They even considered integration with respect to an Euler characteristic of motivic nature and so they introduced the notion of \emph{generalized Poincar\'e series} of a complex curve singularity \cite{CDG10}.

\subsection{} In the spirit of the preceding paragraphs, the author showed in his thesis \cite{M} (see also the joint paper with his advisor Delgado \cite{DM}) that the factors $Z(\mathcal{O}_P, \mathcal{O}_P,T)$ coincide essentially with the generalized Poincar\'e series of Campillo, Delgado and Gusein-Zade, under a suitable specialization for finite fields (see \S \ref{3.8} below). These ideas have also provided some feedback: for instance St\"ohr achieved a deeper insight into the nature of the local zeta functions (see \cite{St3}, and \cite{St4} together with his student J.J. Mira). 

%\S \ref{pa:I1}

\subsection{}  The key ingredient that allows to relate those different formal power series is the \emph{universal zeta function} for a curve singularity defined by Z\'u\~niga and the author in \cite{MZ}: for example, the local zeta functions and Poincar\'e series mentioned above are specializations of this universal zeta function. After some preliminaries, we devote Section \ref{section3} to describe this series. Moreover, we claim that one may establish the local-global behaviour explained in \S \ref{pa:I1} for curves defined over non-finite fields. This conjectural behaviour has already shown some evidences in particular cases; see e.g. the theorem in Section \ref{section4}.

\section{Preliminaries and notations} 

%\textbf{In this section $k$ is a field of characteristic $p\geq0$.???}

\subsection{}\label{2.1} Consider the normalization $\pi: \tilde{X} \to X$ of an algebraic curve $X$ over $k$, and let $\mathcal{O}=\mathcal{O}_P:=\mathcal{O}_{P,X}$ be the local ring of $X$ at $P$. For the sake of simplicity we will assume the ring $\mathcal{O}$ to be complete.
\medskip

It is $\pi^{-1}(P)=\{Q_1,\ldots , Q_d\}$ and so the corresponding local rings $\mathcal{O}_{Q_i}$ are discrete valuation rings of $K$ over $\mathcal{O}$. The value semigroup associated to $\mathcal{O}$ is defined to be
$$
S(\mathcal{O}):=\{\underline{v}(\underline{z}) : \underline{z} \ \mbox{nonzero divisor in } \mathcal{O}\} \subseteq \mathbb{N}^d;
$$
here $\underline{v}(\underline{z})=(v_1(z_1),\ldots , v_d(z_d))$, where each $v_i$ stands for the valuation associated with $\mathcal{O}_{Q_i}$; we write $S$ for this semigroup from now on. Let $c=c(S)$ denote the conductor of $S$, i.e. the smallest element $\underline{v}\in S$ such that $\underline{v}+\mathbb{N}^d\subseteq S$. Moreover, $\mathcal{O}^{\times}$ denotes the group of units of $\mathcal{O}$. Further details here and in the sequel  can be checked in \cite{MZ} and the references therein.

\subsection{} We say that the ring $\mathcal{O}$ is totally rational if all rings $\mathcal{O}_{Q_i}$, for $i=1,\ldots , d$ have $k$ as a residue field.

\subsection{} The integral closure of $\mathcal{O}$ in $K/k$ is $\tilde{\mathcal{O}}=\tilde{\mathcal{O}}_P=\mathcal{O}_{Q_1}\cap \ldots \cap \mathcal{O}_{Q_d}$. We write $\tilde{\mathcal{O}}^{\times}$ for its group of units. The \emph{singularity degree} $\delta_P$ of $\tilde{\mathcal{O}}$ is defined as $\delta_P=\delta:=\dim_k \tilde{\mathcal{O}}/\mathcal{O} <\infty$ (see e.g. \cite[Chapter IV]{Serre}).

\subsection{} For $\underline{n}\in S$ we set
$$
\mathcal{I}_{\underline{n}}:=\left\{  I\subseteq\mathcal{O}\mid I=\underline{z}\mathcal{O}%
\text{, with }\underline{v}(\underline{z})=\underline{n}\right\} ,
$$
and for $m\in\mathbb{N}$,%
\[
\mathcal{I}_{m}:=%
%TCIMACRO{\tbigcup \nolimits_{\substack{\underline{n}\in S\\\left\Vert
%\underline{n}\right\Vert =m}}}%
%BeginExpansion
{\textstyle\bigcup\nolimits_{\substack{\underline{n}\in S\\\left\Vert
\underline{n}\right\Vert =m}}}
%EndExpansion
\mathcal{I}_{\underline{n}},
\]
where $\left\Vert
\underline{n}\right\Vert$ denotes the sum of the components of the vector $\underline{n}=(n_1,\ldots , n_d)\in \mathbb{N}^d$.

\subsection{} In the category $\mathrm{Var}_k$ of $k$-algebraic varieties, we define the Grothendieck ring $K_0(\mathrm{Var}_k)$, which is the ring generated by symbols $[V]$ for $V \in \mathrm{Var}_k$, with the relations $[V]=[W]$ if $V$ is isomorphic to $W$, $[V]=[V\setminus Z]+[Z]$ if $Z$ is closed in $V$, and $[V\times W]=[V][W]$. We write $\mathbb{L}:=[\mathbb{A}^1_{k}]$ for the class of the affine line, and $\mathcal{M}_k:=K_0(\mathrm{Var}_k)[\mathbb{L}^{-1}]$ for the ring obtained by localization with respect to the multiplicative set generated by $\mathbb{L}$.

\subsection{} \label{2.7} It is possible to associate to $\mathcal{I}_{\underline{n}}$~resp.~$\mathcal{I}_m$ well-defined classes in the Grothen\-dieck ring \cite[Section~5]{MZ}; those classes will be denoted by $[\mathcal{I}_{\underline{n}}] $~resp.~$[\mathcal{I}_m]$. This allows to attach to the local ring $\mathcal{O}$ the zeta functions
\[
Z\left(  T_{1},\ldots,T_{d},\mathcal{O}\right)  :=%
%TCIMACRO{\tsum \nolimits_{\underline{n}\in S}}%
%BeginExpansion
{\textstyle\sum\nolimits_{\underline{n}\in S}}
%EndExpansion
\left[  \mathcal{I}_{\underline{n}}\right] \mathbb{L}^{-\left\Vert
\underline{n}\right\Vert }T^{\underline{n}}\in\mathcal{M}_{k} [\![
T_{1},\ldots,T_{d} ]\!]  , \nonumber%
\]
where $T^{\underline{n}}:=T_{1}^{n_{1}} \cdot\ldots\cdot T_{d}^{n_{d}}$, and
$Z\left(  T,\mathcal{O}\right)  :=Z\left(  T,\ldots,T,\mathcal{O}\right) $.

\subsection{} \textbf{Definition.} Consider an algebraic curve $X$ over $k$. If $k$ has characteristic $p\geq0$, then we say that \emph{$k$ is
big enough for $X$} if for every singular point $P$ in $X$ the
following two conditions hold: 1) the ring $\mathcal{O}$ is
totally rational and 2) $\widetilde{\mathcal{O}}^{\times}/\mathcal{O}^{\times} \cong\left(  G_{m}\right)  ^{d-1}%
\times\left(  G_{a}\right)  ^{\delta-d+1}$, with $G_m=(k^{\times},\cdot)$ and $G_a=(k,+)$.
\medskip

Note that  the condition `$k$ is big enough for $X$' is
fulfilled when $p$ is big enough. 

%\begin{proposition}
%\label{lema8} For $\underline{n}\in S$, and $I_0:=\{1, 2, \ldots , r\}$ we have
%$$
%\left[  \mathcal{I}_{\underline{n}}\right]  =\left(
%\mathbb{L}-1\right)  ^{-1}\mathbb{L}^{\left\Vert \underline{n}\right\Vert +1}%
%%TCIMACRO{\tsum \limits_{I\subseteq I_{0}}}%
%%BeginExpansion
%{\textstyle\sum\limits_{I\subseteq I_{0}}}
%%EndExpansion
%\left(  -1\right)  ^{\#(I)}\mathbb{L}^{-\dim_k \mathcal{O}/\{\underline{z} \in \mathcal{O} : \underline{v}(\underline{z}) \geq \underline{n}+\underline
%{1}_{I}  \}}
%$$
%\end{proposition}

\section{The universal zeta function for curve singularities}  \label{section3}

\subsection{} For $k=\mathbb{C}$, we consider a semigroup
$S\subseteq\mathbb{N}^{d}$ such that $S=S\left(  \mathcal{O}\right) $. Moreover, for $\underline{n}\in S$ set 
$$
\mathcal{I}_{\underline{n}}\left(  U\right)  :=\left(  U-1\right)
^{-1}U^{\left\Vert \underline{n}\right\Vert +1}%
{\textstyle\sum\limits_{I\subseteq [d]}}
\left(  -1\right)  ^{\#(I)}U^{-\dim_k \big ( \mathcal{O}/\{\underline{z} \in \mathcal{O} : \underline{v}(\underline{z}) \geqslant \underline{n}+\underline{1}_{I}  \} \big ) },
$$
for an indeterminate $U$, and where $[d]:=\{1,2,\ldots ,d\}$, and $\underline{1}_{I}$ is the $d$-tuple with the components corresponding to the indices in $I$ equal to $1$, and the other components equal to $0$.

The notation $\mathcal{I}_{\underline{n}}\left(  U\right)$ is appropriate, since that expression coincides with $\left[  \mathcal{I}_{\underline{n}}\right]$  when $U$ specializes to $\mathbb{L}$, cf. \cite[Section~5]{MZ}.

\subsection{} 

Let $\underline{c}=(c_1,\ldots, c_d)$ be the conductor of the semigroup  $S$, cf. \S \ref{2.1}. Let $J:=\{1,\ldots , r\}\subseteq [d]$, and let $\underline{m}\in \mathbb{N}^d$ be such that $\underline{c}>\underline{m}$, i.e., $c_i>m_i$ for all $i\in [d]$.  For a fixed $\emptyset \subsetneq J \subsetneq [d]$, set $r_J:=\# J$ and
$$
B_J:=\{\underline{m} \in \mathbb{N}^{r_J}: H_{J,\underline{m}} \neq \emptyset \},
$$
where $H_{J,\underline{m}} :=\{\underline{n} \in S : n_j\geq c_j ~\mbox{if}~j\in J,~\mbox{and}~n_j=m_j~\mbox{otherwise}\}$.
\medskip

\noindent \textbf{Definition}
We define the \emph{universal zeta function} $\mathcal{Z}\left(  T_{1},\ldots,T_{d},U,S\right) $ associa\-ted with $S$ to be
\begin{align*}
{\textstyle\sum\limits_{%
\substack{\underline{n} \in S \\ \underline{0} \leqslant \underline{n} <
\underline{c}} }}
\mathcal{I}_{\underline{n}}\left(  U\right)  U^{-\left\Vert \underline
{n}\right\Vert }T^{\underline{n}}%
%\]%
%\begin{align*}
&+%
{\textstyle\sum\limits_{\emptyset\subsetneq J\subsetneq I_{0}}}
\text{ \ }%
{\textstyle\sum\limits_{\substack{\underline{m}\in B_{J}}}}
\left(  U-1\right)  U^{\left\Vert \underline{c}\right\Vert -\delta
-1}
\mathcal{I}_{\underline{f_J}(\underline{m})}\left(  U\right)
U^{-\left\Vert
\underline{c}\right\Vert -\left\Vert \underline{f_{J}}(\underline{m})\right\Vert }\times \\
\times& \frac{T^{\underline{f_{J}}(\underline{m})}}{%
{\textstyle\prod\limits_{i=1}^{r_{J}}}
\left(  1-U^{-1}T_{i}\right)  }
+\frac{\left(  U-1\right)  ^{d-1}U^{\delta-d+1}U^{-\left\Vert
\underline
{c}\right\Vert } T^{\underline{c}}}{%
{\textstyle\prod\limits_{i=1}^{d}}
\left(  1-U^{-1}T_{i}\right)  },
\end{align*}
where $\underline{f_{J}}(\underline{m})=\left(  c_{1},\ldots,c_{r_{J}},m_{r_{J}+1}%
,\ldots,m_{d}\right)  \in S$, with $m_{i}<c_{i}$, $r_{J}+1\leqslant i\leq d$,
and $1\leqslant r_{J}<d$.

\subsection{}  Observe that this universal zeta function is completely determined by $S$. The adjective \emph{universal} applied to this zeta function will be clear after the follo\-wing paragraphs.

%\subsection{}  Next we justify why to apply the adjective \emph{universal} to this zeta function. First of all, let us observe that the zeta function $Z(T_1,\ldots , T_d,\mathcal{O})$ may be also defined for curves over some other fields than $\mathbb{C}$; in fact, the fields of our interest must contain enough elements. More precisely, we define:\\

\subsection{}  The generalized Poincar\'e series $P_g(T_1,\ldots , T_d)$ of Campillo, Delgado and Gusein-Zade (\cite{CDG10}; see also \cite{CDK}, \cite{DM}) as an integral with respect to an Euler characteristic of motivic nature is very close to the zeta function $Z(T_1,\ldots , T_d, \mathcal{O})$ of \S \ref{2.7}, and therefore to the universal zeta function via the specialization $U=\mathbb{L}$:\\

\noindent \textbf{Proposition.}
\emph{If $S=S(\mathcal{O})$ and $k$ is big enough for $Y$, then}
$$
Z(T_1,\ldots , T_d,\mathcal{O})=\mathbb{L}^{\delta + 1} P_g(T_1,\ldots , T_d) = \mathcal{Z}\left(  T_{1},\ldots,T_{d},U,S\right) |_{U=\mathbb{L}}.
$$

\subsection{}\label{TheoremMONO} \label{TheoremC}  In addition, a certain specialization of the universal zeta function coincides with the zeta function of the monodromy transformation of a reduced plane curve singularity acting on its Milnor fibre, as we briefly explain now.\\

\noindent \textbf{Definition.}
Let $(X,0)\subseteq(\mathbb{C}^{2},0)$ be a reduced plane curve
singularity defined by an equation $f=0$, with
$f\in\mathcal{O}_{(\mathbb{C}^{2},0)}$ reduced. Let
$h_{f}:V_{f}\rightarrow V_{f}$ be the monodromy transformation of
the singularity $f$ acting on its Milnor fiber $V_{f}$. The zeta function of the monodromy $h_{f}$ is defined to be
$$
\varsigma_{f}\left(  T\right):=\prod_{i\geqslant  0} \Big [ \mathrm{det}(\mathrm{id}-T\cdot (h_f)_{\ast} |_{H_i(V_f;\mathbb{C})}) \Big]^{(-1)^{i+1}}.
$$

A result of Campillo, Delgado and Gusein-Zade (\cite[Theorem 1]{CDG1}) allows us to prove:\\

\noindent \textbf{Theorem.}
\emph{Let
$k=\mathbb{C}$. Then for every $\mathcal{O}=\mathcal{O}_{\left(  \mathbb{C}%
^{2},0\right)  }/\left(  f\right) $, with $f\in\mathcal{O}_{\left(
\mathbb{C}^{2},0\right)  }$ reduced, and for every $S=S\left(  \mathcal{O}%
\right)  $, we have}
\[
\varsigma_{f}\left(  T\right)  =\mathcal{Z}\left(  T_{1},\ldots,T_{d}%
,U,S\right)  \mid_{%
\begin{array}
[c]{l}%
T_{1}=\ldots=T_{d}=T\\
U=1
\end{array}
.}%
\]

\subsection{}  In \cite{Z} Z\'u\~niga introduced a Dirichlet series $Z(\mathrm{Ca}%
(X),T)$ associated to the effective Cartier divisors on an algebraic curve $X$
defined over $k=\mathbb{F}_{q}$, which admits
an Euler product of the form%
\[
Z(\mathrm{Ca}(X),T)=%
%TCIMACRO{\tprod \limits_{P\in X}}%
%BeginExpansion
{\textstyle\prod\limits_{P\in X}}
%EndExpansion
Z_{\mathrm{Ca}(X)}(T,q,\mathcal{O}_{P,X}),
\]
with $Z_{\mathrm{Ca}(X)}(T,q,\mathcal{O}_{P,X}):=
{\textstyle\sum\limits_{I\subseteq \mathcal{O}_{P,X}}}
T^{\dim_{k}\left(\mathcal{O}_{P,X}/I\right)  }$, 
where $I$ runs through all the principal ideals of $\mathcal{O}_{P,X}$.  In addition, $Z_{\mathrm{Ca}%
(X)}(T,q,\mathcal{O}_{P,X})=Z(T,\mathcal{O}_{P,X})$, cf.~\S \ref{2.7}.
\medskip

Observe that this zeta function is nothing but the zeta function $Z(\mathcal{O},\mathcal{O},T)$ appearing as a local factor in the St\"ohr zeta function, cf.~ Section \ref{intro}.

\subsection{} \textbf{Remark.}\label{3.8}
In the category of $\mathbb{F}_{q}$-algebraic varieties, $\left[
\cdot\right]  $ specializes to the counting rational points additive invariant
$\#\left(  \cdot\right)  $. We write $\#\left(  Z\left(  T_{1},\ldots
,T_{d},\mathcal{O}\right)  \right)  $ for the rational function obtained by
specializing $\left[  \cdot\right]  $ to $\#\left(  \cdot\right)  $. From a
computational point of view, $\#\left(  Z\left(  T_{1},\ldots,T_{d}%
,\mathcal{O}\right)  \right)  $ is obtained from $Z\left(  T_{1},\ldots
,T_{d},\mathcal{O}\right)  $ by replacing $\mathbb{L}$ by $q$.

\subsection{}\label{TheoremD} \textbf{Theorem.}
{\it Let $\ k=\mathbb{F}_{q}$ and let $\mathcal{Z}\left(
T_{1},\ldots,T_{d},U,S\right)  $ be the universal zeta function for $S$. Moreover, let
$X$ be an algebraic curve defined over $k$, and let $\mathcal{O}_{P,X}$ be the
(complete) local ring of $X$ at a singular point $P$. Assume that $k$
is big enough for $X$ and that $S=S\big(\mathcal{O}_{P,X}\big)  $.\\
\\
\noindent(1) For any $\mathcal{O}=\mathcal{O}_{\left(  \mathbb{C}%
^{2},0\right)  }/\left(  f\right)  $, with $f\in\mathcal{O}_{\left(
\mathbb{C}^{2},0\right)  }$ reduced, and $S=S(\mathcal{O})$,%
\[
Z_{\mathrm{Ca}(X)}\big(  q^{-1}T,q,\mathcal{O}_{P,X}\big)    =\#\left(
Z\big(  T_{1},\ldots,T_{d},\mathcal{O}_{P,X}\big)  \right) 
\]
\[
=\mathcal{Z}\big(  T_{1},\ldots,T_{d},U,S\big)  \mid_{%
\begin{array}
[c]{l}%
{\small T_{1}=\ldots=T_{d}=T}\\
{\small U=q}
\end{array}
.}%
\]
In particular $Z_{\mathrm{Ca}(X)}\big(  q^{-1}T,q,\mathcal{O}_{P,X}\big)  $
depends only on $S$. Moreover, if $X$\ is plane, then
$Z_{\mathrm{Ca}(X)}\big(  q^{-1}T,q,\mathcal{O}_{P,X}\big)  $ is a complete
invariant of the equisingularity class of $\mathcal{O}_{P,X}$.

\noindent(2) For any $\mathcal{O}=\mathcal{O}_{\left(  \mathbb{C}%
^{2},0\right)  }/\left(  f\right)  $, with $f\in\mathcal{O}_{\left(
\mathbb{C}^{2},0\right)  }$, it holds that
\[
Z_{\mathrm{Ca}(X)}\big(  q^{-1}T,q,\mathcal{O}_{P,X}\big)  \mid
_{q \rightarrow 1}=\varsigma_{f}( T).
\]
}

\section{Some connections between local and global zeta functions}\label{section4}

\subsection{} For a smooth algebraic variety $Y$ defined over a field $k$, M. Kapranov defined a zeta function as the formal power series in an indeterminate $u$
$$
\zeta_{\mathrm{mot},Y} (u)=\sum_{n=0}^{\infty} [Y^{(n)}] u^n \in K_0(\mathrm{Var}_k)[\![u]\!],
$$
where $Y^{(n)}$ stands for the
$n$-fold symmetric product of $Y$ (cf. \cite[\S 1]{K}). (For instance, if $k=\mathbb{F}_q$, then one obtains the usual Hasse-Weil zeta function of  $Y$, cf.~\S \ref{3.8}). When $Y$ is a curve, Baldassarri, Deninger and Naumann introduced in \cite{BDN} a two-variable version of the Kapranov zeta function, namely
$$
Z_{\mathrm{mot},Y}(t,u)=\sum_{n,\nu \geqslant 0} [\mathrm{Pic}_{\nu}^{n}] \frac{u^{\nu}-1}{u-1} t^n \in K_0(\mathrm{Var}_k)[\![u,t]\!],
$$
where the algebraic $k$-scheme $\mathrm{Pic}_{\nu}^{n}=\mathrm{Pic}_{\geqslant \nu}^{n} \setminus \mathrm{Pic}_{\geqslant \nu+1}^{n}$ (with $\mathrm{Pic}_{\geqslant \nu}^{n} $ being the closed subvariety---in the Picard variety of degree $n$ line bundles on $Y$---of line bundles $\mathcal{L}$ with $h^0(\mathcal{L}) \geqslant \nu$) defines a class in $K_0(\mathrm{Var}_k)$.
\medskip

\subsection{} The connections between the universal zeta function and the motivic zeta functions of Kapranov and Baldassarri-Deninger-Naumann are being currently investigated by A. Melle, W. Z\'u\~niga and the author; we believe that the zeta functions discussed in the previous sections are factors of motivic zeta functions of Baldassarri-Deninger-Naumann type for singular curves (as mentioned before, this is known when $k=\mathbb{F}_q$). In order to give some evidence supporting this belief, this note will be finished by stating the relation between local and global zeta functions in a particular situation.
\medskip

The context will be the one of a \emph{modulus}: Following Serre \cite{Serre}, let $k$ be an algebraically closed field, and let $C$ be an irreducible, non-singular, complete algebraic curve defined over $k$. If $F$ is a
finite subset of $C$, a modulus $\mathfrak{m}$ supported on $F$ is defined to be
the assignment of an integer $n_{P} >0$ for each point $P \in F$; this is
sometimes identified with the effective divisor $\sum_{P} n_{P} P$. 
\medskip

\subsection{} It is possible to attach a curve to $\mathfrak{m}$ starting from $C$, essentially by ``placing'' the points in $F$ all together into one (see again \cite{Serre}). The resulting singular curve $C_{\mathfrak{m}}$ has then this point as its only singularity. It holds the following\\

\noindent \textbf{Theorem.}
{\it Let $C_{\mathfrak{m}}$ be a curve arising from a modulus $\mathfrak{m}$ supported on a finite set of points of a curve $C$ as above, and let $P$ be its only singular point.  Furthermore, let $\pi:\widetilde{C_{\mathfrak{m}}} \to C_{\mathfrak{m}}$ be the normalization morphism. Then
$$
Z_{\mathrm{mot},C_{\mathfrak{m}}}(\mathbb{L}^{-1}T,\mathbb{L})= Z_{\mathrm{mot},\widetilde{C_{\mathfrak{m}}}}(\mathbb{L}^{-1}T,\mathbb{L}) \prod_{i=1}^{\sharp(\pi^{-1}(P))} (1-\mathbb{L}^{-1} T)
\cdot Z (T, \mathcal{O}_{P}).
$$
}
The proof of this statement will appear in a forthcoming paper.

\end{document}